\documentclass[reqno]{amsart}
\usepackage{iftex}
\ifPDFTeX
    \usepackage[utf8]{inputenc}
    \usepackage{lmodern}
\else
    \usepackage[no-math]{fontspec}
    \defaultfontfeatures{Ligatures=TeX}
\fi
\usepackage{amsmath, amsthm, amssymb}
\usepackage{aliascnt}
\usepackage{mathrsfs}

\usepackage[dvipsnames]{xcolor}
\definecolor{linkBlue}{rgb}{0,.25,.55}

\usepackage[psdextra]{hyperref}
\hypersetup{
    colorlinks=true,
    linkcolor=linkBlue,
    filecolor=magenta,
    urlcolor=linkBlue,
    citecolor=linkBlue}
\usepackage[nameinlink,noabbrev]{cleveref}

\usepackage{mdframed}
\usepackage{upgreek}
\usepackage{multicol}
\usepackage{float}
\usepackage{adjustbox}

\newcommand{\Hom}{\text{Hom}}

\usepackage{graphicx}
\graphicspath{ {./images/} }
\usepackage[english]{babel}
\usepackage[justification=centering]{subfig}
\usepackage{spverbatim}
\usepackage{mathtools}
\usepackage[shortlabels]{enumitem}
\usepackage{pgfplots}
\pgfplotsset{compat=1.18}
\usepackage{tikz-cd}
\usepackage{stmaryrd}

\usepackage[rightcaption]{sidecap}

\newcommand{\iprod}{\mathbin{\raisebox{.2ex}{$\lrcorner$}}}

\theoremstyle{plain}
\newtheorem{theorem}{Theorem}[section]
\newcommand{\namedtheoremname}{}
\newtheorem*{namedtheoreminner}{\namedtheoremname}
\newenvironment{namedtheorem}[1]{
    \renewcommand{\namedtheoremname}{#1}
    \begin{namedtheoreminner}
}{
    \end{namedtheoreminner}
}
\newaliascnt{proposition}{theorem}
\newtheorem{proposition}[proposition]{Proposition}
\aliascntresetthe{proposition}
\newaliascnt{lemma}{theorem}
\newtheorem{lemma}[lemma]{Lemma}
\aliascntresetthe{lemma}
\newaliascnt{corollary}{theorem}
\newtheorem{corollary}[corollary]{Corollary}
\aliascntresetthe{corollary}

\theoremstyle{definition}
\newaliascnt{definition}{theorem}
\newtheorem{definition}[definition]{Definition}
\aliascntresetthe{definition}
\newaliascnt{assumption}{theorem}
\newtheorem{assumption}[assumption]{Assumption}
\aliascntresetthe{assumption}
\crefname{assumption}{assumption}{assumptions}
\Crefname{assumption}{Assumption}{Assumptions}
\newaliascnt{example}{theorem}
\newtheorem{example}[example]{Example}
\aliascntresetthe{example}

\theoremstyle{remark}
\newaliascnt{remark}{theorem}
\newtheorem{remark}[remark]{Remark}
\aliascntresetthe{remark}

\newcommand \Q {\mathbb Q}
\newcommand \C {\mathbb C}

\newcommand \Z {\mathbb Z}

\let\tilde\widetilde

\numberwithin{equation}{section}

\newcommand{\de}{\mathrm{d}}

\newcommand{\prim}{\mathrm{prim}}
\newcommand{\van}{\mathrm{van}}

\DeclareMathOperator{\tr}{tr}

\DeclareMathOperator{\GL}{GL}
\DeclareMathOperator{\SL}{SL}
\DeclareMathOperator{\PSL}{PSL}

\usepackage{csquotes}
\usepackage[
    backend=biber,
    style=alphabetic,
    sorting=nyt,
    maxalphanames=3,
    minalphanames=1
]{biblatex}

\DeclareSourcemap{
  \maps[datatype=bibtex]{
    \map{
      \step[fieldsource=doi, final]
      \step[fieldsource=url, match=\regexp{^https?://(dx\.)?doi\.org/}, final]
      \step[fieldset=url, null]
      \step[fieldset=urldate, null]
    }
  }
}

\AtBeginBibliography{\setlength{\emergencystretch}{2em}}
\usepackage{bbm}

\usepackage[margin={1.25in,1.5in}]{geometry}
\usepackage{thmtools}
\usepackage{cancel}

\DeclareMathOperator{\Res}{Res}
\definecolor{ssNakano}{HTML}{3A6591}
\definecolor{ssDiagonal}{HTML}{1D742B}
\definecolor{ssAdjacent}{HTML}{912B65}

\newcommand\Gr{\mathrm{Gr}}

\newcommand{\acts}{\mathbin{{\text{\raisebox{.75em}{\rotatebox{-90}{$\circlearrowright$}}}}}}

\title{Griffiths--Dwork Reduction For Smooth Hypersurfaces in Grassmannians}

\author{Andrew Huchala}
\address{University of Oregon, USA}
\email{ahuchala@uoregon.edu}

\author{Arya Yae}
\email{ayae@uoregon.edu}



\begin{document}

\begin{abstract}
We give an explicit Griffiths--Dwork pole-order reduction algorithm for a smooth hypersurface $Z$ in a Grassmannian $\Gr(k,n)$. Under two arithmetic conditions on $k,n$, and the degree of $Z$, the residue construction describes the Hodge filtration on the middle cohomology of the complement in terms of the Jacobian ring. We explicitly characterize the required cohomological vanishings. As applications, we compute group representations on the vanishing cohomology of two symmetric hyperplane sections of $\Gr(3,10)$ and deduce maximal Picard rank for their associated Debarre--Voisin fourfolds. We also compute a Picard--Fuchs equation for a pencil of symmetric Gushel--Mukai fivefolds.
\end{abstract}

\maketitle

\tableofcontents

\section{Introduction and Motivation}
\label{sec: intro}
\subsection{Overview}
The Hodge structure of a smooth hypersurface $Z\subset \mathbb P^n$ may be understood by Griffiths' theory of residues \cite{Gri69}. To each associated graded Hodge piece of the primitive (or equivalently, vanishing) cohomology of $Z$ we associate a particular graded part of a quotient ring $R/J_f$, given in terms of the defining equation $f$ for $Z$. 

The ring $R/J_f$ is called the \emph{Jacobian ring} due to its particularly simple presentation as a quotient by the partial derivatives of $f$:
\[J_f := \left(\frac{\partial}{\partial x_0}f,\cdots, \frac{\partial}{\partial x_n}f\right)
\]
where $R=\mathbb C[x_0,\cdots,x_n]$ is the homogeneous coordinate ring of $\mathbb P^n$. The residue-theoretic construction has computational advantages which we outline here.

Every class in $H^n(\mathbb P^n\setminus Z)$ is represented by a rational form $\frac{g\Omega}{f^p}$ with $\Omega\in H^0(K_{\mathbb P^n}(n+1))$, and one can relate the pole order $p$ to the corresponding piece of the Hodge filtration. If $g$ lies in the \emph{Jacobian ideal} $J_f$, say $g = \sum h_i \frac{\partial f}{\partial x_i}$, then for $p>1$,
\[
  \dfrac{g\Omega}{f^p}\equiv \dfrac{1}{p-1}\dfrac{\sum \frac{\partial h_i}{\partial x_i}\Omega}{f^{p-1}}\qquad\text{modulo exact forms.}
\]
A polynomial with nonzero image in $R/J_f$ cannot be reduced, and its class gives the corresponding associated Hodge graded piece via the residue map. This is \emph{Griffiths--Dwork reduction}.  By choosing a graded complement to $J_f$ in $R$, the reduction can be iterated to decompose the class $[g\Omega/f^p]$ as a sum of representatives of decreasing pole order.

Saito \cite[Definition~3.1.2]{Sai87} extended the definition of the Jacobian ring beyond $\mathbb P^n$ to a smooth hypersurface $Z=V(f)$ in a compact irreducible Hermitian symmetric space $X$.  Let $R = \bigoplus_{a\geq 0}H^0(X,\mathcal O_X(a))$ denote the homogeneous coordinate ring of $X$. Then the Jacobian ring of $Z$ is defined as $R/J_f$, where
\[
  J_f=(f)+(v(f) : v\in H^0(X,T_X)).
\]
We consider the case where $X$ is a Grassmannian $\Gr(k,n)$.

\begin{namedtheorem}{\Cref{thm: main theorem}}
    Let $Z=V(f)$ be a smooth hypersurface of degree $d$ in $X=\Gr(k,n)$, where $2\leq k\leq n-2$, and put $N=k(n-k)$.  Assume that
    \begin{itemize}
        \item if $d\mid n$, then $\gcd(k,n/d)=1$, and
        \item if $k\in\{2,n-2\}$, then $2d\nmid(n+1)$.
    \end{itemize}
    Then for every $p\geq1$ there is an isomorphism
    \[
        (R/J_f)_{pd-n}\;\cong\;H^{N-p,\,p-1}_{\mathrm{van}}(Z)\oplus
        H^{N-p+1,\,p-1}_{\mathrm{prim}}(X).
    \]
\end{namedtheorem}

The second summand is a contribution from the ambient Grassmannian, absent in the projective case, which is only potentially nonzero when $p = N/2+1$ (and hence $N$ is even).

Fatighenti and Mongardi \cite{Fat21} prove the existence of the above isomorphisms, given certain cohomological vanishing hypotheses on $X$.  These vanishings are guaranteed when $d \geq n-1$ \cite[Theorem~3.2 and Corollary in \S~3.5]{Sno86}, giving the isomorphism in those cases.  \Cref{thm: main theorem} requires extending to $d\leq n-2$.  We apply the classification of the relevant cohomology groups obtained by Gossow and the first author \cite[Theorems 4.6 and 4.8]{Gos26}. In \Cref{prop: vanishings}, we characterize these cohomological vanishings using the arithmetic hypotheses in the theorem.  This makes the theory usable in low degree, where many interesting examples were previously excluded by the uniform bound $d \geq n-1$.  For example, families of Debarre--Voisin twentyfolds occur at $d=1$, and families of Gushel--Mukai fivefolds occur at $d=2$.

In addition to characterizing the vanishing conditions, our construction of the isomorphism works on the level of meromorphic forms, leading to an explicit Griffiths--Dwork reduction algorithm: \Cref{thm:pole order reduction} states that for $p>1$ and a certain global section $\Omega$ of the twisted canonical bundle $K_X(n)$, a meromorphic form $g\Omega/f^p$ of suitable degree is cohomologous to some $g'\Omega/f^{p-1}$ if $g\in J_f$.  The converse holds under the same arithmetic assumptions of \Cref{thm: main theorem}.  A similar reduction identity appears in Huang, Lian, Yau and Yu \cite[(2.6)--(2.7)]{Hua18}, who also identify the Hodge filtration of the cohomology of the complement in terms of the Jacobian ring under cohomological vanishing assumptions \cite[Theorem 2.4]{Hua18}.

\bigskip

In \Cref{sec: GIT}, we use the GIT construction of $X=\Gr(k,n)$ to explicitly construct a global section $\Omega$ of the twisted canonical bundle $K_X(n)$ and prove that it is nowhere vanishing and invariant under the natural action $\SL(n,\C)\acts X$. 

In \Cref{sec: hypersurfaces and residues}, we prove a residue theorem for smooth ample hypersurfaces in smooth projective varieties under two cohomological vanishing hypotheses. The proof compares truncated logarithmic de Rham complexes with complexes of meromorphic forms of bounded pole order, then applies the hypercohomology spectral sequence. We characterize the vanishing hypotheses arithmetically for Grassmannians and derive the differential identity used for pole-order reduction.

In \Cref{sec: pole order reduction}, we prove the main theorem and give an explicit Griffiths--Dwork reduction algorithm. We also show how to perform the reduction within individual isotypic components of cohomology under finite groups of linear symmetries. 

\Cref{sec: Applications} demonstrates a few calculations utilizing the results of the previous sections. For a $G$-invariant hypersurface $Z$, the primitive cohomology of $Z$ decomposes into irreducible representations of $G$, which our algorithm may actually compute. We may use this in conjunction with \Cref{thm: main theorem} to compute the decomposition of Hodge summands into irreducible $G$-representations. For example, we present a new proof of the maximal Picard rank of a certain Debarre--Voisin fourfold with $\PSL_2(11)$ symmetry considered by Song \cite{Son21}. 

Song computes the integral Picard lattice of this fourfold, and also recalls Mongardi's result that a symplectic automorphism of order $11$ forces maximal Picard rank \cite[Theorem 5.10]{Son21}. Here we recover the rank and the relevant representations from character theory applied to the Jacobian ring of the associated twentyfold, a hyperplane section of $\Gr(3,10)$. We also determine the relevant Hodge representations for another Debarre--Voisin fourfold with symmetry group $A_7$. Both groups occur in the classification of symplectic actions by H\"ohn and Mason \cite{HM19}; Wawak \cite[Example 3.6]{Waw23} also studies the $\PSL_2(11)$ Debarre--Voisin example. Here the Jacobian ring gives a direct calculation of the Hodge representations and hence of maximal Picard rank.

We also apply the reduction algorithm to a pencil of $A_5$-invariant quadric sections of $\Gr(2,5)$. By passing to the subspace fixed by $A_4$, we obtain a variation of rank two and compute its Picard--Fuchs equation. This example illustrates how symmetry and pole-order reduction may be used together to make a period calculation manageable.

\subsection{Acknowledgements}
The authors wish to thank Nicolas Addington and Fern Gossow for many productive discussions. This material is based upon work supported by the National Science Foundation under Grant No. DMS-2424556 while the first author was in residence at the Institute for Computational and Experimental Research in Mathematics in Providence, RI, during the Computations on K3 Surfaces and Related Varieties program. Both authors were also partially supported by the University of Oregon.

\emph{AI Disclosure.} This paper grew out of a series of discussions around the first author's thesis. The ideas and proof techniques are our own. Large language models were used for proofreading and coding assistance.
\section{A Construction of a Global Section of \texorpdfstring{$K_X(n)$}{K(n)} via GIT}
\label{sec: GIT}
In this section, we produce a global section of the twisted canonical bundle $K_X(n)$ and show it is invariant under the natural action $\SL(n,\C)\acts X$.  We work with GIT construction of the Grassmannian $X=\Gr(k,n)$ as the quotient of the stable locus $M^\mathrm{st}$ in $M:=\Hom(\mathbb C^k,\mathbb C^n)$ by the right action of $\GL(k,\mathbb{C})$ on $M$, using the character $\chi(g)=\det(g)$.

Let $(x_{11},x_{12},\ldots,x_{nk})$ be the standard coordinates on $M$, and let $\tilde\Omega\in H^0(M,\Omega_M^{nk})$ be the algebraic top form
\[\tilde\Omega = \de x_{11}\wedge\de x_{12}\wedge\ldots\wedge \de x_{nk}.\]
Any $A\in\mathfrak{gl}(k,\C)$ generates a vertical vector field $\zeta_A$ on $M$.  Since the right action of $e^{tA}$ on $\tilde\Omega$ scales by $\det(e^{tA})^n$, differentiating yields
\begin{equation}\label{eq: lie derivative omega tilde}
    \mathcal L_{\zeta_A}\tilde\Omega=n\tr(A)\tilde\Omega.
\end{equation}
On the other hand, $\SL(n,\mathbb C)$ acts by left multiplication on $M$ and an element $A\in\mathfrak{sl}(n,\C)$ generates a vector field $\zeta_A$ on $M$.  An analogous argument to the above shows that
\begin{equation}\label{eq: lie derivative is trace}
    \mathcal L_{\zeta_A}\tilde\Omega=k\tr(A)\tilde\Omega=0.
\end{equation}

The vertical vectors generated by the elementary matrices $f_{ij}\in\mathfrak{gl}(k,\mathbb C)$ can be written as
\[F_{ij} = \sum_{\ell = 1}^n x_{\ell i}\,\dfrac{\partial}{\partial x_{\ell j}}.\]
Define the contracted form
\[\Omega := F_{11} \iprod F_{12}\iprod\ldots\iprod F_{kk}\iprod \tilde\Omega.
\]
\begin{lemma}
    \label{lem: omega global section}
    $\Omega$ is a semi-invariant form on $M$ with weight $n$ with respect to the $\GL(k,\mathbb C)$-action, and descends to a global section of $K_X(n)$.
\end{lemma}
\begin{proof}
    Write $R_g(x)=xg$ for the right action of $g\in\GL(k,\mathbb C)$.  A form $f$ on the stable locus $M^{\text{st}}$ is semi-invariant with weight $n$ if and only if $R_g^*f = (\det g)^n f$ for all $g\in\GL(k,\mathbb C)$.
    We first prove the equivalent linearized statement: for any $B\in\mathfrak{gl}(k,\C)$ we have
    \begin{equation}\label{eq: infinitesimal semiinvariance property}
        \mathcal L_{\zeta_B}\Omega = n\,\text{Tr} (B)\,\Omega.
    \end{equation}
    Expanding $\mathcal L_{F_{ij}}\Omega$ yields
    \[\mathcal L_{F_{ij}}\Omega = \mathcal L_{F_{ij}}\left(F_{11} \iprod F_{12}\iprod\ldots\iprod F_{kk}\iprod \tilde\Omega\right).
    \]
    Recalling the identity $[\mathcal L_X,\iota_Y]\alpha = \iota_{[X,Y]}\alpha$ (see e.g. \cite[\S1.3]{Kob63}),
    \begin{align*}
        \mathcal L_{F_{ij}}&\left(F_{11} \iprod F_{12}\iprod\cdots\iprod F_{kk}\iprod \tilde\Omega\right)\\
        &=[F_{ij},F_{11}] \iprod F_{12}\iprod\cdots\iprod F_{kk}\iprod \tilde\Omega + F_{11}\iprod  \mathcal L_{F_{ij}}\left(F_{12}\iprod\cdots\iprod F_{kk}\iprod \tilde\Omega\right).
    \end{align*}
    Repeating this to commute the Lie derivative past the contractions,
    \begin{equation}\label{eq: Lie derivative of Omega}
        \mathcal L_{\zeta_B}\Omega=\left(\sum_{i,j}T_{ij}\right)+F_{11}\iprod\cdots\iprod F_{kk}\iprod\mathcal L_{\zeta_B}\tilde\Omega
    \end{equation}
    where
    \[T_{ij}:=F_{11}\iprod\cdots\iprod[\zeta_B,F_{ij}]\iprod\cdots\iprod F_{kk}\iprod\tilde\Omega.\]
    The last term in \eqref{eq: Lie derivative of Omega} is the desired $n\tr(B)\Omega$ by \eqref{eq: lie derivative omega tilde}.  It remains to show $T_{ij}=-T_{ji}$ for all $i,j$, hence the sum $\sum_{i,j}T_{ij}=0$.  Since this condition is linear in $B$, it suffices to assume $B$ is an elementary matrix $f_{rs}$, in which case $\zeta_B=F_{rs}$.  Then
    \[[\zeta_B,F_{ij}] = \delta_{si}F_{rj} - \delta_{rj}F_{is}.\]
    Since contraction is alternating, $T_{ij}$ is zero unless the coefficient of $F_{ij}$ in this expansion is nonzero, which requires $r=s=i$ or $r=s=j$.  Therefore we may assume $r=s$, in which case
    \begin{align*}
        T_{ij}
        &=F_{11}\iprod\cdots\iprod(\delta_{si}F_{sj}-\delta_{sj}F_{is})\iprod\cdots\iprod F_{kk}\iprod\tilde\Omega\\
        &=F_{11}\iprod\cdots\iprod(\delta_{si}F_{ij}-\delta_{sj}F_{ij})\iprod\cdots\iprod F_{kk}\iprod\tilde\Omega.
    \end{align*}
    Now it is clear that $T_{ij}=-T_{ji}$.  This proves \eqref{eq: infinitesimal semiinvariance property}.  Since $\GL(k,\mathbb C)$ is connected, integrating gives $R_g^*\Omega=(\det g)^n\Omega$.  On the stable locus $M^{\mathrm{st}}$, the fields $F_{ij}$ span the tangent space to the $\GL(k,\C)$-orbit, so $\zeta_A\iprod\Omega=0$ by alternation for all $A\in\mathfrak{gl}(k,\C)$.  Thus $\Omega$ descends to a global section of $K_X(n)$.
\end{proof}

At a point $x\in M^\mathrm{st}$, extend $F_{11},\dots,F_{kk}$ to a basis $b_{11},\dots,b_{nk}$ for $T_xM\cong M$.  Then
\[
    \Omega(b_{k+1,1},\dots,b_{nk})=\tilde\Omega(b_{11},\dots,b_{nk})\ne0.
\]
Thus the descended form $\Omega$ on $\Gr(k,n)$ is nowhere vanishing.

The vector fields $\zeta_A$ arising from the linearization of the left action of $\SL(n,\C)\acts M$ descend to vector fields on $X$ (which we also denote by $\zeta_A$) because the actions of $\SL(n,\C)$ and $\GL(k,\C)$ commute.
\begin{proposition}
    \label{prop: lie derivative Eij}
    For any $A\in\mathfrak{sl}(n,\C)$,
    \[\mathcal L_{\zeta_A}(\Omega) = 0.\]
\end{proposition}
\begin{proof}
    Since the actions of $\SL(n,\C)$ and $\GL(k,\C)$ on $M$ commute, we have $[\zeta_A,F_{ij}]=0$ for all $i,j$.  The identity $[\mathcal L_{\zeta_A},\iota_{F_{ij}}]=\iota_{[\zeta_A,F_{ij}]}=0$ therefore lets us commute $\mathcal L_{\zeta_A}$ past each contraction.
    \begin{align*}
        \mathcal L_{\zeta_A}(\Omega) = \mathcal L_{\zeta_A}\left(F_{11} \iprod F_{12}\iprod\ldots\iprod F_{kk}\iprod \tilde\Omega\right)
        &=F_{11} \iprod F_{12}\iprod\ldots\iprod F_{kk}\iprod \mathcal L_{\zeta_A}\tilde\Omega=0,
    \end{align*}
    where the last equality follows from \eqref{eq: lie derivative is trace} and the fact that $\tr(A)=0$.
\end{proof}

\section{Hypersurfaces and Residues}
\label{sec: hypersurfaces and residues}

In this section, we express the vanishing cohomology of a smooth ample hypersurface in terms of meromorphic forms on its complement.  Under two cohomological vanishing assumptions, \Cref{thm: Griffiths} identifies forms of a fixed pole order, modulo lower-order poles and exact forms, with the corresponding Hodge components of the vanishing cohomology of the hypersurface and the primitive cohomology of the ambient variety.  We first prove this residue theorem for a general smooth ample hypersurface.  For Grassmannians, \Cref{prop: vanishings} gives numerical criteria for the two assumptions.  We end the section with a useful reduction identity for exact forms.

\subsection{Residues, the Hodge Filtration, and Pole Order}
\label{sec:hypercohom log de rham}

Let $Z\subset X$ be a smooth ample hypersurface in a smooth projective variety of dimension $N$.  We consider the following two vanishing assumptions.

\begin{assumption}
    \label{asm:ij-N}
    $H^i(X,\Omega_X^{N-i}((p-i)Z))=0$ for all $0<i<p\leq N$.
\end{assumption}

\begin{assumption}
    \label{asm:ij-N-minus-1}
    $H^i(X,\Omega_X^{N-i-1}((p-i-1)Z))=0$ for all $0<i<p-1<N$.
\end{assumption}

Set $U=X\setminus Z$, and let $\ell:Z\hookrightarrow X$ and $j:U\hookrightarrow X$ denote the inclusions.  The long exact sequence of the pair $(X,U)$ gives
\[
    \cdots\longrightarrow H^m(X,U,\mathbb Q)\longrightarrow H^m(X,\mathbb Q)\longrightarrow H^m(U,\mathbb Q)\longrightarrow H^{m+1}(X,U,\mathbb Q)\longrightarrow\cdots
\]
The Thom isomorphism identifies $H^{m+1}(X,U,\mathbb Q)$ with $H^{m-1}(Z,\mathbb Q)(-1)$.  Under this identification, the map to $H^{m+1}(X,\mathbb Q)$ is the Gysin morphism $\ell_*$, and the connecting homomorphism from $H^m(U,\mathbb Q)$ is the residue map.  Since the residue removes one differential, it maps $F^pH^m(U)$ into $F^{p-1}H^{m-1}(Z)$.  The Tate twist $(-1)$ accounts for this shift by setting
\[
    F^p\bigl(H^{m-1}(Z)(-1)\bigr)=F^{p-1}H^{m-1}(Z),
\]
so that the residue map preserves the grading.  In degree $N$ we obtain
\[
    \begin{aligned}
        H^{N-2}(Z,\mathbb Q)(-1)&\xrightarrow{\ell_*}H^N(X,\mathbb Q)\xrightarrow{j^*}H^N(U,\mathbb Q)\xrightarrow{\Res}H^{N-1}(Z,\mathbb Q)(-1)\xrightarrow{\ell_*}H^{N+1}(X,\mathbb Q).
    \end{aligned}
\]
Define
\[
    H^{N-1}_{\mathrm{van}}(Z,\mathbb Q) = \ker\bigl(\ell_*:H^{N-1}(Z,\mathbb Q)\longrightarrow H^{N+1}(X,\mathbb Q)\bigr).
\]
Since $Z$ is ample, restriction to $Z$ identifies $H^{N-2}(X,\mathbb Q)$ with $H^{N-2}(Z,\mathbb Q)$ by the Lefschetz hyperplane theorem.  Under this identification, the Gysin map is given by the cup product with $[Z]$, and thus the Lefschetz decomposition identifies its cokernel with $H^N_{\mathrm{prim}}(X,\mathbb Q)$.  Thus the Gysin sequence gives a short exact sequence
\begin{equation}
    \label{eq:vanishingCohom}
    0\longrightarrow H^N_{\mathrm{prim}}(X,\mathbb Q) \longrightarrow H^N(U,\mathbb Q) \xrightarrow{\Res} H^{N-1}_{\mathrm{van}}(Z,\mathbb Q)(-1) \longrightarrow 0.
\end{equation}
These are morphisms of mixed Hodge structures, so taking $\operatorname{Gr}_F^{N-p+1}$ preserves this short exact sequence.  For $1\leq p\leq N$, rewriting the outer terms we have
\begin{equation}
    \label{eq: filtered HNU}
    \begin{split}
        0\longrightarrow H^{N-p+1,p-1}_{\mathrm{prim}}(X)\longrightarrow \operatorname{Gr}_F^{N-p+1}H^N(U,\mathbb C)\longrightarrow H^{N-p,p-1}_{\mathrm{van}}(Z)\longrightarrow 0.
    \end{split}
\end{equation}
Next, we'll describe $\Gr_F^{N-p+1}H^N(U,\C)$ in terms of meromorphic forms and their pole orders.

We introduce the logarithmic de Rham complex to compute the cohomology of $U$, together with its Hodge filtration. Recall that $\Omega_X^k(\log Z)$ consists of meromorphic $k$-forms $\alpha$ for which both $\alpha$ and $\de\alpha$ have at most simple poles along $Z$. In local coordinates with $Z=\{z_1=0\}$,
\[
    \Omega_X^1(\log Z) = \mathcal O_X\frac{\de z_1}{z_1} \oplus \bigoplus_{i=2}^N\mathcal O_X\,\de z_i,
    \qquad \Omega_X^k(\log Z)=\bigwedge^k\Omega_X^1(\log Z).
\]
These sheaves are locally free, and the exterior derivative of a logarithmic form is again logarithmic.

The logarithmic de Rham complex computes the cohomology of the complement (see, e.g. \cite[Corollary 8.19]{Voi01}):
\[
    \mathbb H^N(X,\Omega_X^\bullet(\log Z))\cong H^N(U,\mathbb C).
\]
The Hodge filtration is induced by truncation of the logarithmic de Rham complex.  Since $X$ is smooth and projective and $Z$ is smooth, Deligne's degeneration theorem \cite[Corollaire 3.2.13(ii)]{Del71} implies that the logarithmic Hodge--de Rham spectral sequence degenerates at $E_1$.  This shows that the natural map
\[\mathbb H^N(X,\Omega_X^{\geq k}(\log Z))\longrightarrow H^N(U,\mathbb C)
\]
is injective, with image $F^kH^N(U,\mathbb C)$.  Here $\Omega_X^{\geq k}(\log Z)$ is obtained by replacing $\Omega_X^j(\log Z)$ by zero for $j<k$.  At $k=0$ this is the full logarithmic complex, so $F^0H^N(U,\mathbb C)=H^N(U,\mathbb C)$.

To relate the Hodge filtration to pole order, for $0\leq k\leq N$ consider the complex $\mathcal C_k^\bullet$ defined by
\[
    \mathcal C_k^\bullet=\left[ \Omega_X^k(Z) \xrightarrow{\de}\Omega_X^{k+1}(2Z)\xrightarrow{\de}\cdots\xrightarrow{\de}\Omega_X^N((N-k+1)Z)\right],
\]
supported in degrees $k,\ldots,N$.

\begin{lemma}
    \label{lem: pole order resolution}
    For $0\leq k\leq N$, the natural inclusion
    \[
        \Omega_X^{\geq k}(\log Z)\hookrightarrow\mathcal C_k^\bullet
    \]
    is a quasi-isomorphism.
\end{lemma}

\begin{proof}
    Suppose first that $1\leq k\leq N$.  Since nothing maps into degree $k$, the cohomology sheaf of $\mathcal C_k^\bullet$ in that degree is the sheaf of closed $k$-forms with at most a simple pole along $Z$.  Exactness in degree $r$, for $k<r\leq N$, amounts to showing that every closed local section $\alpha$ of $\Omega_X^r((r-k+1)Z)$ can locally be written as $\de\beta$ with $\beta$ a section of $\Omega_X^{r-1}((r-k)Z)$.  Both $r$ and $r-k+1$ are at least two, so the existence of $\beta$ follows from \cite[Corollary 6.7]{Voi02}.

    Every closed form with at most a simple pole is logarithmic, and the logarithmic de Rham complex is locally exact in degrees at least two by \cite[Proposition 8.18]{Voi01}.  Thus $\Omega_X^{\geq k}(\log Z)$ and $\mathcal C_k^\bullet$ resolve the same sheaf of closed $k$-forms in degree $k$, so the inclusion is a quasi-isomorphism.

    For the $k=0$ case, since $Z$ is smooth, Deligne's filtered comparison theorem identifies the full logarithmic complex $\Omega_X^\bullet(\log Z)$ with $\mathcal C_0^\bullet$ up to quasi-isomorphism at filtration level zero \cite[\S3]{Del71}.
\end{proof}

The above lemma, together with the identification of $F^kH^N(U,\mathbb C)$ with the hypercohomology of $\Omega_X^{\geq k}(\log Z)$, shows
\begin{equation}
    \label{eq: F equals hyper}
    F^{N-p+1}H^N(U,\mathbb C) \cong \mathbb H^N(X,\mathcal C_{N-p+1}^\bullet),\qquad 1\leq p\leq N+1.
\end{equation}

For $1\leq p\leq N$, we now apply \Cref{asm:ij-N,asm:ij-N-minus-1} to compute the hypercohomology of $\mathcal C_{N-p+1}^\bullet$ using global meromorphic forms.  In the hypercohomology spectral sequence, \Cref{asm:ij-N} eliminates the higher-cohomology terms of total degree $N$, and \Cref{asm:ij-N-minus-1} eliminates the incoming differentials $d_r$ for $r\geq2$ into the remaining term $E_r^{N,0}$.

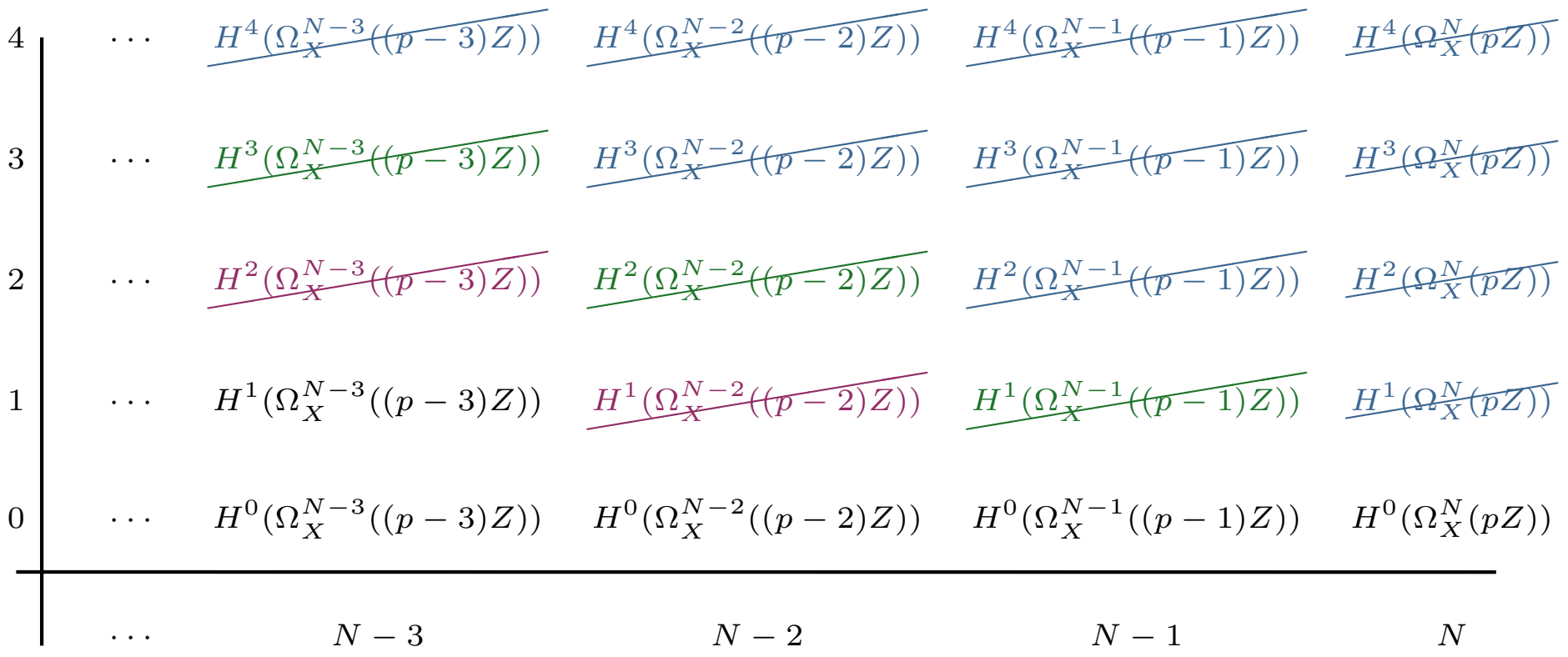
\begin{figure}[H]
    \centering
    \begin{tikzpicture}
        \matrix (m) [matrix of math nodes,
            nodes in empty cells,
            nodes={minimum width=5ex, minimum height=5ex, outer sep=-5pt, font=\scriptsize},
            column sep=1ex,row sep=1ex]{
            4 & \cdots & {\color{ssNakano}\cancel{H^4(\Omega_X^{N-3}((p-3)Z))} }& {\color{ssNakano}\cancel{H^4(\Omega_X^{N-2}((p-2)Z))}} & {\color{ssNakano}\cancel{H^4(\Omega_X^{N-1}((p-1)Z))}} & {\color{ssNakano}\cancel{H^4(\Omega_X^{N}(pZ))}} \\
            3 & \cdots & {\color{ssDiagonal}\cancel{H^3(\Omega_X^{N-3}((p-3)Z))}} & {\color{ssNakano}\cancel{H^3(\Omega_X^{N-2}((p-2)Z))}} & {\color{ssNakano}\cancel{H^3(\Omega_X^{N-1}((p-1)Z))}} & {\color{ssNakano}\cancel{H^3(\Omega_X^{N}(pZ))}} \\
            2 & \cdots & {\color{ssAdjacent}\cancel{H^2(\Omega_X^{N-3}((p-3)Z))}} & {\color{ssDiagonal}\cancel{H^2(\Omega_X^{N-2}((p-2)Z))}} & {\color{ssNakano}\cancel{H^2(\Omega_X^{N-1}((p-1)Z))}} & {\color{ssNakano}\cancel{H^2(\Omega_X^{N}(pZ))}} \\
            1 &\cdots & H^1(\Omega_X^{N-3}((p-3)Z)) & {\color{ssAdjacent}\cancel{H^1(\Omega_X^{N-2}((p-2)Z))}} & {\color{ssDiagonal}\cancel{H^1(\Omega_X^{N-1}((p-1)Z))}} & {\color{ssNakano}\cancel{H^1(\Omega_X^{N}(pZ))}} \\
            0 &\cdots & H^0(\Omega_X^{N-3}((p-3)Z)) & H^0(\Omega_X^{N-2}((p-2)Z)) & H^0(\Omega_X^{N-1}((p-1)Z)) & H^0(\Omega_X^{N}(pZ)) \\
            \quad &\cdots & N-3 & N-2 & N-1 & N \\
        };

        \draw[thick] (m-1-1.east) -- (m-6-1.east);
        \draw[thick] ([yshift=-1.5ex]m-5-1.south) -- ([xshift=1em,yshift=-1.5ex]m-5-6.south);

    \end{tikzpicture}
    \caption{\sloppy The cohomological vanishings on the $E_1$ page of the hypercohomology spectral sequence.  Terms {\color{ssNakano}above the antidiagonal} vanish by Nakano.  Terms {\color{ssDiagonal}on the antidiagonal with $q>0$} vanish by \Cref{asm:ij-N}.  Terms {\color{ssAdjacent}immediately to the left of the antidiagonal with $q>0$} vanish by \Cref{asm:ij-N-minus-1}.}
    \label{fig: spectral sequence 5x5 full}
\end{figure}

\begin{proposition}
    \label{prop: pole order Hodge filtration}
    Under \Cref{asm:ij-N,asm:ij-N-minus-1}, the map
    \begin{equation}
        \label{eq: mapAlpha}
        \alpha_p:H^0(X,K_X(pZ))\longrightarrow H^N(U,\mathbb C),\qquad \omega\longmapsto[\omega|_U],
    \end{equation}
    induces an isomorphism for $2\leq p\leq N$
    \begin{equation}
        \label{eq: hypercohomology degen}
        \frac{H^0(X,K_X(pZ))}{\de H^0(X,\Omega_X^{N-1}((p-1)Z))}\xrightarrow{\sim}F^{N-p+1}H^N(U,\mathbb C).
    \end{equation}
    For $p=1$, it induces an isomorphism $H^0(X,K_X(Z))\cong F^NH^N(U,\mathbb C)$.
\end{proposition}

\begin{proof}
    For $p\geq2$, consider the spectral sequence
    \[
        E_1^{a,b} = H^b(X,\Omega_X^a((a-N+p)Z)) \Longrightarrow \mathbb H^{a+b}(X,\mathcal C_{N-p+1}^\bullet),
    \]
    where $N-p+1\leq a\leq N$. \Cref{asm:ij-N} forces every term of total degree $N$ except $E_1^{N,0}$ to be zero, and
    \[
        E_2^{N,0} = \frac{H^0(X,K_X(pZ))}{\de H^0(X,\Omega_X^{N-1}((p-1)Z))}.
    \]
    There are no outgoing differentials from this term. For $r\geq2$, the source of an incoming differential is a
    subquotient of
    \[
        E_1^{N-r,r-1} = H^{r-1}(X,\Omega_X^{N-r}((p-r)Z)).
    \]
    This vanishes by \Cref{asm:ij-N-minus-1} when $r<p$, and lies outside the complex when $r\geq p$. Thus $E_2^{N,0}=E_\infty^{N,0}$ is the only surviving term of total degree $N$. The resulting edge isomorphism is induced by $\alpha_p$, thus proving the claim by \eqref{eq: F equals hyper}.

    For $p=1$, the complex has the single term $K_X(Z)$ in degree $N$.
\end{proof}

Taking successive quotients of the Hodge filtration gives the following residue theorem.

\begin{restatable}{theorem}{thmGriffiths}
    \label{thm: Griffiths}
    Let $Z\subset X$ be a smooth ample hypersurface in a smooth projective variety of dimension $N$.  Assume the cohomological vanishings of \Cref{asm:ij-N,asm:ij-N-minus-1}.  For $2\leq p\leq N$, there is an isomorphism of complex vector spaces
    \begin{equation}\label{eq: isom to van and prim}
        \frac{H^0(X,K_X(pZ))}{H^0(X,K_X((p-1)Z))+\de H^0(X,\Omega_X^{N-1}((p-1)Z))}
        \cong H^{N-p,p-1}_{\mathrm{van}}(Z)
        \oplus H^{N-p+1,p-1}_{\mathrm{prim}}(X).
    \end{equation}
    For $p=1$,
    \[
        H^0(X,K_X(Z))
        \cong
        H^{N-1,0}_{\mathrm{van}}(Z)
        \oplus
        H^{N,0}_{\mathrm{prim}}(X).
    \]
    Here $\de$ denotes the exterior derivative. 
\end{restatable}

\begin{proof}
    For $p\geq2$, the natural inclusions
    \[
        \Omega_X^{\geq N-p+2}(\log Z) \hookrightarrow \Omega_X^{\geq N-p+1}(\log Z),
        \qquad \mathcal C_{N-p+2}^\bullet \hookrightarrow \mathcal C_{N-p+1}^\bullet
    \]
    commute with the comparison maps of \Cref{lem: pole order resolution}. Hence the isomorphisms in \Cref{prop: pole order Hodge filtration} identify $F^{N-p+2}H^N(U,\mathbb C)$ with the image of $H^0(X,K_X((p-1)Z))$. Taking the quotient gives
    \begin{equation}\label{eq: graded hodge pieces as quotient}
        \operatorname{Gr}_F^{N-p+1}H^N(U,\mathbb C) \cong\frac{H^0(X,K_X(pZ))}{H^0(X,K_X((p-1)Z))+\de H^0(X,\Omega_X^{N-1}((p-1)Z))}.
    \end{equation}
    The result for $2\leq p\leq N$ follows from choosing a splitting of the exact sequence \eqref{eq: filtered HNU}.  For the case $p=1$, the statement follows from the observation $F^NH^N(U,\mathbb C)\cong H^0(X,K_X(Z))$ and $F^{N+1}H^N(U,\mathbb C)=0$.
\end{proof}

\subsection{Numerical Criteria for Grassmannians}

We now specialize to Grassmannians.  Let $X=\Gr(k,n)$, with $2\leq k\leq n-2$, and let $Z=V(f)$ be a smooth hypersurface of degree $d$.  Write $N=k(n-k)$.  The vanishing hypotheses admit the following numerical characterization.

\begin{proposition}
    \label{prop: vanishings}
    \Cref{asm:ij-N,asm:ij-N-minus-1} are characterized as follows:
    \begin{enumerate}
        \item \label{item: one}
        \Cref{asm:ij-N} fails if and only if $d\mid n$ and $\gcd(k,n/d)>1$.
        \item \label{item: two}
        \Cref{asm:ij-N-minus-1} fails if and only if $k\in\{2,n-2\}$ and $2d\mid(n+1)$.
    \end{enumerate}
\end{proposition}
\begin{proof}
    First we prove (\ref{item: one}). By \cite[Theorem 4.6]{Gos26}, the positive twists $t$ for which $H^i(X,\Omega_X^{N-i}(t))\neq0$ for some $i>0$ are precisely
\[
t=\frac{n}{s},\qquad s>1,\qquad s\mid\gcd(k,n).
\]
Such a twist is divisible by $d$ precisely when $sd\mid n$. Thus one exists if and only if $d\mid n$ and $\gcd(k,n/d)>1$. The corresponding cohomological degree satisfies $i<N/2$, while $t\leq n/2\leq N/2$, so the pole order $p=i+t/d$ lies in the required range.

Now we prove (\ref{item: two}).
By \cite[Theorem 4.8]{Gos26}, the only nonzero groups $H^i(X,\Omega_X^{N-i-1}(t))$ with $i,t>0$ occur when $k\in\{2,n-2\}$ and $n$ is odd, with
\[
i=\frac{n-3}{2},\qquad t=\frac{n+1}{2}.
\]
This twist is divisible by $d$ precisely when $2d\mid(n+1)$. The corresponding pole order satisfies
\[
i+1<p=i+1+\frac{t}{d}\leq i+1+t=n\leq N.
\]
\end{proof}

\begin{example}
    A smooth hyperplane section of $\Gr(3,10)$ satisfies both hypotheses: $\gcd(3,10)=1$ and $3\notin\{2,8\}$. In contrast, a smooth hyperplane section of $\Gr(2,5)$ fails the second hypothesis, since $2\mid6$. A smooth quadric section of $\Gr(2,5)$ satisfies both hypotheses, since $2\nmid5$ and $4\nmid6$.
\end{example}

\begin{remark}\label{rem: all positive twists}
    The same calculation gives a slightly stronger statement which will be useful for high pole orders.  When \Cref{asm:ij-N,asm:ij-N-minus-1} hold,
    \[
        H^i(X,\Omega_X^{N-i}(md))=H^i(X,\Omega_X^{N-i-1}(md))=0
    \]
    for every $m>0$ and $i>0$.
\end{remark}

\subsection{Exact Forms on Grassmannians}

Next, we will describe the exact forms in \Cref{thm: Griffiths} for a hypersurface $Z=V(f)$ in $X=\Gr(k,n)$.  Write $N=k(n-k)$ and $d=\deg f$ as before.  Recall that $K_X\cong\mathcal O_X(-n)$ and $T_X\cong S^\vee\otimes Q$, where
\begin{equation}
    \label{eq: tautological}
    0\longrightarrow S\longrightarrow V\otimes\mathcal O_X \longrightarrow Q\longrightarrow0, \qquad V=\mathbb C^n
\end{equation}
is the tautological exact sequence, with $S$ the rank-$k$ tautological subbundle and $Q$ the rank-$(n-k)$ quotient bundle.

Let $\{A_i\}$ be a basis of $\mathfrak{sl}(n,\mathbb C)$, and let $\zeta_{A_i}$ denote the corresponding fundamental vector fields.  To express global forms in terms of these vector fields, we will need two vanishing results.  See \cite[Corollaries 2.12 and 2.13]{HuchalaThesis} for the proofs.

\begin{lemma}
    \label{lem: Qdual vanishing}
    For any line bundle $L$ on $X=\Gr(k,n)$ with $2\leq k\leq n-1$,
    \[
        H^1(X,Q^\vee\otimes V\otimes K_X\otimes L)=0.
    \]
\end{lemma}

\begin{lemma}
    \label{lem: EndS vanishing}
    For any line bundle $L$ on $X=\Gr(k,n)$ with $1\leq k\leq n-2$,
    \[
        H^1(X,\operatorname{End}(S)\otimes K_X\otimes L)=0.
    \]
\end{lemma}

For $2\leq k\leq n-2$, both vanishings apply and give the following description.

\begin{proposition}
    \label{cor: varphi}
    For $p\geq1$, every global section of $\Omega_X^{N-1}(pZ)$ can be written as
    \[
        \sum_i\frac{P_i(\zeta_{A_i}\iprod\Omega)}{f^p}, \qquad \deg P_i=pd-n,
    \]
    where $\Omega$ is the generator of $H^0(X,K_X(n))$ constructed in \Cref{sec: GIT}.
\end{proposition}

\begin{proof}
    Dualizing the tautological sequence \eqref{eq: tautological} and tensoring with $K_X(pZ)\otimes V$ gives
    \begin{equation}
        \label{eq: Qdual}
        0\longrightarrow K_X(pZ)\otimes V\otimes Q^\vee
        \longrightarrow K_X(pZ)\otimes\operatorname{End}(V)
        \xrightarrow{\varphi_1}K_X(pZ)\otimes V\otimes S^\vee
        \longrightarrow0.
    \end{equation}
    On the other hand, tensoring \eqref{eq: tautological} with $K_X(pZ)\otimes S^\vee$ gives
    \begin{equation}
        \label{eq: endS}
        0\longrightarrow K_X(pZ)\otimes\operatorname{End}(S)
        \longrightarrow K_X(pZ)\otimes V\otimes S^\vee
        \xrightarrow{\varphi_2}K_X(pZ)\otimes T_X
        \longrightarrow0.
    \end{equation}
    Applying \Cref{lem: Qdual vanishing,lem: EndS vanishing} with $L=\mathcal O_X(pZ)$ shows that the first cohomology of each left-hand term vanishes.  The long exact cohomology sequences for \eqref{eq: Qdual} and \eqref{eq: endS} therefore show that $\varphi_1$ and $\varphi_2$ induce surjections on global sections.  Define the composition
    \[
        \varphi:H^0(X,K_X(pZ))\otimes\operatorname{End}(V)
        \twoheadrightarrow H^0(X,K_X(pZ)\otimes T_X).
    \]
    Both $\varphi_1$ and $\varphi_2$ are defined using the identity on the tensor factor $K_X(pZ)$.  The composition on the remaining tensor factors takes an endomorphism $A$ to the fundamental vector field $\zeta_A$.  Thus $\varphi(\omega\otimes A)=\omega\otimes\zeta_A$ for all $\omega\otimes A$ in the domain.

    Scalar matrices act trivially on $X$, so restricting $\varphi$ to $H^0(X,K_X(pZ))\otimes\mathfrak{sl}(n,\mathbb C)$ remains surjective.
    By the contraction isomorphism
    \[
        H^0(X,K_X(pZ)\otimes T_X)\xrightarrow{\sim}H^0(X,\Omega_X^{N-1}(pZ)),
    \]
    every global section $\eta$ of $\Omega_X^{N-1}(pZ)$ corresponds to an element $\xi\in H^0(X,K_X(pZ)\otimes T_X)$.  An arbitrary element in the domain of the restricted map $\varphi$ has the form
    \[\sum_iP_i\Omega/f^p\otimes A_i\]
    where $P_i\in H^0(X,\mathcal O_X(pd-n))$ and $\{A_i\}$ is our basis for $\mathfrak{sl}(n,\C)$.  Surjectivity of $\varphi$ therefore lets us write
    \[
        \xi=\varphi\left(\sum_i\frac{P_i\Omega}{f^p}\otimes A_i\right)
        =\sum_i\frac{P_i\Omega}{f^p}\otimes\zeta_{A_i}.
    \]
    Applying the contraction isomorphism again gives $\eta=\sum_i P_i(\zeta_{A_i}\iprod\Omega)/f^p$, as desired.
\end{proof}

The following identity produces the pole-order reduction used
in the next section.

\begin{lemma}
    \label{lem: GD reduction}
    Let $p\geq2$, let $P$ be homogeneous of degree $(p-1)d-n$, and let $A\in\mathfrak{sl}(n,\mathbb C)$. Then
    \[
        \de\left(\frac{P(\zeta_A\iprod\Omega)}{f^{p-1}}\right)
        =\frac{\zeta_A(P)\Omega}{f^{p-1}}-(p-1)\frac{P\zeta_A(f)\Omega}{f^p}.
    \]
\end{lemma}

\begin{proof}
    Set $\omega=P\Omega/f^{p-1}$. Since $\omega$ is a meromorphic top form, $\de\omega=0$. Cartan's formula therefore gives
    \[
        \de(\zeta_A\iprod\omega)=\mathcal L_{\zeta_A}\omega.
    \]
    The result follows from the Leibniz rule and the invariance of $\Omega$ established in \Cref{prop: lie derivative Eij}.
\end{proof}

\section{Griffiths--Dwork Reduction for Grassmannians}
\label{sec: pole order reduction}
In this section, we prove our main theorem relating the Jacobian ring of a Grassmannian hypersurface to the vanishing cohomology of that hypersurface.  We then give a Griffiths--Dwork reduction algorithm, along with a version that preserves isotypic cohomology components under the action of a finite group $G\subset\GL(n,\C)$.

\subsection{Residue Isomorphism}
\label{sec: Residue Isomorphism}
Consider the action of $\mathfrak{sl}(n,\mathbb C)\cong H^0(X,T_X)$ on the homogeneous coordinate ring $R$ of $X=\Gr(k,n)$ given by differentiation along the vector fields $\zeta_A$. This action preserves the grading. Let $f\in R_d$ define a smooth hypersurface $Z$, and set $U=X\setminus Z$.

\begin{definition}\label{def: Jacobian ideal}\label{def: Jacobian ring}
    The \emph{Jacobian ideal} $J_f$ in $R$ is the homogeneous ideal generated by $f$ and all $\zeta_A(f)$ for $A\in\mathfrak{sl}(n,\mathbb C)$. The \emph{Jacobian ring} is the graded ring $R/J_f$.
\end{definition}

We can now prove our first main result.
\begin{theorem}\label{thm: main theorem}
    Let $X=\Gr(k,n)$ with $2\leq k\leq n-2$, let $N=k(n-k)$, and let $Z=V(f)\subset X$ be a smooth hypersurface of degree $d$. Assume that
    \begin{enumerate}
        \item if $d\mid n$, then $\gcd(k,n/d)=1$; and
        \item if $k\in\{2,n-2\}$, then $2d\nmid(n+1)$.
    \end{enumerate}
    For every $p\geq1$, the map $g\mapsto[g\Omega/f^p]$ induces an isomorphism
    \[
        (R/J_f)_{pd-n}\xrightarrow{\sim}\operatorname{Gr}_F^{N-p+1}H^N(U,\mathbb C).
    \]
    In light of \eqref{eq: filtered HNU}, this gives the vector space decomposition
    \[
        (R/J_f)_{pd-n}\cong H_{\mathrm{van}}^{N-p,p-1}(Z)\oplus H_{\mathrm{prim}}^{N-p+1,p-1}(X).
    \]
    For $p>N$, both sides are zero.
\end{theorem}
\begin{proof}
    We first consider $1\leq p\leq N$. By \Cref{prop: vanishings}, the hypotheses imply the cohomological vanishings needed in \Cref{thm: Griffiths}. For $p\geq2$, let
    \begin{equation}\label{eq: Kp}
        \mathcal K_p=H^0(X,K_X((p-1)Z))+\de H^0(X,\Omega_X^{N-1}((p-1)Z))
    \end{equation}
    be the kernel of the map induced by $\alpha_p$ to $\operatorname{Gr}_F^{N-p+1}H^N(U,\mathbb C)$, with $\alpha_p$ as defined in \Cref{prop: pole order Hodge filtration}.  Equation \eqref{eq: graded hodge pieces as quotient} gives $H^0(X,K_X(pZ))/\mathcal K_p\cong\operatorname{Gr}_F^{N-p+1}H^N(U,\mathbb C)$. Consider the isomorphism
    \[
        \psi_p:R_{pd-n}\longrightarrow H^0(X,K_X(pZ)),\qquad g\longmapsto\frac{g\Omega}{f^p}.
    \]
    We must show that $\psi_p((J_f)_{pd-n})=\mathcal K_p$.

    Recall that $(J_f)_{pd-n}$ is spanned by terms $gf$ and $g\zeta_A(f)$ with $g\in R_{(p-1)d-n}$. We have
    \[
        \psi_p(gf)=\frac{g\Omega}{f^{p-1}}\in\mathcal K_p.
    \]
    Meanwhile, the identity in \Cref{lem: GD reduction} gives
    \[
        \psi_p(g\zeta_A(f))=\frac{1}{p-1}\left(\frac{\zeta_A(g)\Omega}{f^{p-1}}-\de\left(\frac{g(\zeta_A\iprod\Omega)}{f^{p-1}}\right)\right)\in\mathcal K_p.
    \]
    Thus $\psi_p((J_f)_{pd-n})\subseteq\mathcal K_p$.

    For the reverse containment, take an arbitrary $\delta+\de\gamma\in\mathcal K_p$. By \Cref{cor: varphi}, given a basis $\{A_i\}$ of $\mathfrak{sl}(n,\mathbb C)$, we can write
    \[
        \delta=\frac{g\Omega}{f^{p-1}},\qquad \gamma=\sum_i\frac{P^i(\zeta_{A_i}\iprod\Omega)}{f^{p-1}},\qquad g,P^i\in R_{(p-1)d-n}.
    \]
    Then $\delta=\psi_p(gf)$, and
    \[
        \de\gamma=\psi_p\left(\sum_i\bigl(f\zeta_{A_i}(P^i)-(p-1)P^i\zeta_{A_i}(f)\bigr)\right).
    \]
    Both $gf$ and $\sum_i\bigl(f\zeta_{A_i}(P^i)-(p-1)P^i\zeta_{A_i}(f)\bigr)$ lie in $J_f$, proving the reverse containment.

    When $p=1$, all generators of $J_f$ have degree $d$, so $(J_f)_{d-n}=0$.  Thus the isomorphism follows from the $p=1$ case of \Cref{thm: Griffiths}.  When $p=N+1$, \Cref{rem: all positive twists} gives the vanishings needed in the spectral sequence of \Cref{prop: pole order Hodge filtration}, so the same quotient calculation gives $(R/J_f)_{(N+1)d-n}\cong\operatorname{Gr}_F^0 H^N(U)$.  By $E_1$-degeneration of the logarithmic Hodge--de Rham spectral sequence in \Cref{sec:hypercohom log de rham}, $\operatorname{Gr}_F^0 H^N(U)\cong H^N(X,\mathcal O_X)=0$, since $X$ is a Grassmannian.  All higher degrees of $R/J_f$ vanish since the ring is generated in degree one.
\end{proof}

\subsection{Reduction of Pole Order}
We now turn the differential identity from \Cref{lem: GD reduction} into an explicit reduction step.  We need not assume the arithmetic hypotheses of \Cref{thm: main theorem}.

\begin{theorem}[Griffiths--Dwork reduction]\label{thm:pole order reduction}
    For $p>1$, suppose $g\in(J_f)_{pd-n}$. Write
    \[
        g=cf+\sum_i h^i\zeta_{A_i}(f),\qquad c,h^i\in R_{(p-1)d-n},
    \]
    where $\{A_i\}$ is a basis of $\mathfrak{sl}(n,\mathbb C)$. Then, modulo exact forms on $U$,
    \begin{equation}\label{eq: beta reduced}
        \frac{g\Omega}{f^p}\equiv\frac{\left(c+\frac1{p-1}\sum_i\zeta_{A_i}(h^i)\right)\Omega}{f^{p-1}}.
    \end{equation}
\end{theorem}
\begin{proof}
    Set $\gamma=\sum_i h^i(\zeta_{A_i}\iprod\Omega)/f^{p-1}$. By \Cref{lem: GD reduction}, the difference between the two displayed forms is $-\de\gamma/(p-1)$.
\end{proof}

Suppose $\beta=g\Omega/f^p$ is a meromorphic top form on $X$. Restricting to $U$ gives its cohomology class $[\beta]$. Under the arithmetic hypotheses, \Cref{thm: main theorem} shows that its image in $\operatorname{Gr}_F^{N-p+1}H^N(U)$ vanishes precisely when $g\in J_f$. In that case, \Cref{thm:pole order reduction} explicitly produces a cohomologous form with a smaller pole. For $p>N$, every such form $g\Omega/f^p$ has $g\in J_f$, so this reduction always applies.

\subsection{The Reduction Algorithm}
We normalize the numerators by factorials so that successive reductions have no denominators depending on the pole order.

\begin{definition}\label{def: normalized pole class}
    For $p\geq1$ and $g\in R_{pd-n}$, set
    \[
        [g]_p=\left[\frac{(p-1)!\,g\Omega}{f^p}\right]\in H^N(U,\C).
    \]
    We say that $g$ represents the class $[g]_p$ at pole order $p$.
\end{definition}

In this notation, \Cref{thm:pole order reduction} gives
\begin{equation}\label{eq: factorial reduction}
    \left[cf+\sum_i h^i\zeta_{A_i}(f)\right]_p
    =\left[(p-1)c+\sum_i\zeta_{A_i}(h^i)\right]_{p-1},
\end{equation}
for $p>1$ and $c,h^i\in R_{(p-1)d-n}$.

Assume the hypotheses of \Cref{thm: main theorem}.  Choose in each degree a complement
\[
    R_m=(J_f)_m\oplus W_m,
\]
setting $W_m=0$ for $m<0$.  Fix a basis $\{A_i\}$ of $\mathfrak{sl}(n,\C)$.  Given $g_p\in R_{pd-n}$ with $p\geq1$, reduce $[g_p]_p$ as follows.
\begin{enumerate}[itemsep=0.5em]
    \item For $j=p,p-1,\ldots,2$:
    \begin{enumerate}[label=(\roman*),itemsep=0.5em]
        \item Decompose $g_j=a_j+b_j$, with $a_j\in W_{jd-n}$ and $b_j\in(J_f)_{jd-n}$.
        \item Write $b_j=c_jf+\sum_i h_j^i\zeta_{A_i}(f)$, with $c_j,h_j^i\in R_{(j-1)d-n}$.
        \item Set $g_{j-1}=(j-1)c_j+\sum_i\zeta_{A_i}(h_j^i)$.
    \end{enumerate}
    \item Set $a_1=g_1$ and return $(a_1,\ldots,a_p)$.
\end{enumerate}
The last step is valid because $(J_f)_{d-n}=0$.  If $p>N$, then $a_j=0$ for $N<j\leq p$ by \Cref{thm: main theorem}.

\begin{proposition}\label{prop: reduction algorithm}
    The output of the reduction algorithm satisfies
    \[
        [g_p]_p=\sum_{j=1}^p[a_j]_j.
    \]
    For fixed complements $W_m$, the tuple $(a_1,\ldots,a_p)$ depends only on the class $[g_p]_p$, and not on the expressions chosen in step (ii).
\end{proposition}
\begin{proof}
    Each step is possible by the definition of $J_f$.  Equation \eqref{eq: factorial reduction} gives
    \[
        [g_j]_j=[a_j]_j+[g_{j-1}]_{j-1},
    \]
    and iterating proves the decomposition.

    For uniqueness, suppose that $\sum_{j=1}^p[a_j]_j=0$, with $a_j\in W_{jd-n}$.  If some $a_j$ is nonzero, choose the largest such $j$.  Passing to $\operatorname{Gr}_F^{N-j+1}H^N(U)$ kills all terms with smaller index.  By \Cref{thm: main theorem}, the remaining term can vanish only if $a_j\in J_f$, contradicting the choice of complement.  Thus every $a_j$ is zero.
\end{proof}

Let $W = \bigoplus_m W_m$. The resulting classes $[a_j]_j$ depend on the graded splitting $R=J_f\oplus W$ and need not be of pure Hodge type.  The projection of $[a_p]_p$ to the associated graded is independent of the splitting, while the lower components generally are not. Nevertheless, once the representatives are fixed, the algorithm gives coordinates for every class represented by a rational top form.

\begin{corollary}\label{cor: invariant reduction}
    Assume the hypotheses of \Cref{thm: main theorem}, and let $G\subset\mathrm{GL}(n,\C)$ be a finite subgroup fixing $f$ and $\Omega$.  For an invariant $g_p\in R_{pd-n}^G$, the reduction algorithm can be performed entirely in the invariant subspaces, using complements
    \[
        R_m^G=(J_f)_m^G\oplus W_m.
    \]
    Every numerator $g_j$ and every term $a_j,b_j$ can be chosen invariant.  The resulting tuple $(a_1,\ldots,a_p)$ depends only on $[g_p]_p$ and these complements.
\end{corollary}
\begin{proof}
    Recall from \eqref{eq: factorial reduction}
    \begin{equation*}
        \left[cf+\sum_i h^i\zeta_{A_i}(f)\right]_p
        =\left[(p-1)c+\sum_i\zeta_{A_i}(h^i)\right]_{p-1}
    \end{equation*}
    Since $G$ fixes $f$ and $\Omega$, averaging over $G$ results in invariant representatives on both sides.  If $g\in(J_f)_{pd-n}$ was already invariant, then averaging preserves the left-hand side and produces a lower pole order reduction whose representative belongs to $R_{(p-1)d-n}^G$.  The group action preserves the space of vector fields induced by $\mathfrak{sl}(n,\mathbb C)$, so $J_f$ is $G$-stable.  Thus, all the steps in the algorithm can be carried out without leaving the invariant subspace.
\end{proof}

Now suppose a finite subgroup $G\subset\GL(n,\C)$ acts on $f$ and $\Omega$ by scalars, i.e.\ $f$ and $\Omega$ are respectively $\xi_f$- and $\xi_\Omega$-semi-invariant for characters $\xi_f,\xi_\Omega:G\to\C^\times$.  Then $G$ preserves the hypersurface $Z=\{f=0\}$, and the middle cohomology $H^N(U,\C)$ decomposes into irreducible $G$-representations.  For an irreducible representation $V$ of $G$ with character $\chi$, we write $R_m[V]$ for the sum of all subrepresentations of $R_m$ isomorphic to $V$, and similarly for $(J_f)_m[V]$ and $H^N(U,\C)[V]$.  We call these the \emph{$V$-isotypic components}.  Replace the averaging operation with the isotypic projection
\[
    P_V(g)=\frac{\dim V}{|G|}\sum_{\sigma\in G}\chi(\sigma^{-1})(\sigma\cdot g),
\]
using the same notation $P_V$ for the isotypic projection $H^N(U,\C)\to H^N(U,\C)[V]$.  With $V$ fixed, define $V_j$ to be the twisted representation with character $\chi_j=\chi\xi_f^j\xi_\Omega^{-1}$.  Then for $g\in R_{jd-n}$, we have
\[
    \begin{aligned}
        P_V([g]_j)
        &=\frac{\dim V}{|G|}\sum_{\sigma\in G}\chi(\sigma^{-1})
        \left[\frac{(j-1)!\,(\sigma\cdot g)(\sigma\cdot\Omega)}{(\sigma\cdot f)^j}\right]\\
        &=\frac{\dim V}{|G|}\sum_{\sigma\in G}\chi(\sigma^{-1})\xi_f(\sigma)^{-j}\xi_\Omega(\sigma)
        \left[\frac{(j-1)!\,(\sigma\cdot g)\Omega}{f^j}\right]\\
        &=\left[\frac{\dim V}{|G|}\sum_{\sigma\in G}
        \chi_j(\sigma^{-1})(\sigma\cdot g)\right]_j\\
        &=[P_{V_j}(g)]_j.
    \end{aligned}
\]
By applying $P_V$, every class in $H^N(U,\C)[V]$ represented at pole order $p$ can instead be written as $[g_p]_p$ with $g_p\in R_{pd-n}[V_p]$.

\begin{corollary}\label{cor: semi-invariant reduction}
    With the notation above, assume the hypotheses of \Cref{thm: main theorem}.  Let $\beta=[g_p]_p$ for $g_p\in R_{pd-n}[V_p]$.  Choose complements to $(J_f)_{jd-n}[V_j]$ in $R_{jd-n}[V_j]$ for $j=1,\dots,p$.  Then the reduction algorithm can be performed with each $g_j,a_j,b_j\in R_{jd-n}[V_j]$.
\end{corollary}
\begin{proof}
    Use the same algorithm, projecting at each reduction step as follows.  At the $j$th step, suppose $b_j\in(J_f)_{jd-n}[V_j]$ and the ordinary reduction gives $[b_j]_j=[g_{j-1}']_{j-1}$ for some $g_{j-1}'\in R_{(j-1)d-n}$.  Applying $P_V$ to both sides gives
    \[
        [b_j]_j=P_V([b_j]_j)=P_V([g_{j-1}']_{j-1})=[P_{V_{j-1}}(g_{j-1}')]_{j-1},
    \]
    so we set $g_{j-1}=P_{V_{j-1}}(g_{j-1}')$.  The ideal $J_f$ is still $G$-invariant, so iterating gives the result.
\end{proof}

For smooth hypersurfaces in projective space, the same argument applies to classical Griffiths--Dwork reduction \cite{Gri69}.  See \cite{DKSSVW}.

\section{Applications}
\label{sec: Applications}

The following applications use the accompanying Sage scripts \texttt{debarre\_voisin.sage} and \texttt{pencil\_\allowbreak calculations.sage} to compute group actions on cohomology and a Picard--Fuchs equation, respectively.

\subsection{Background}
Beauville \cite{Beau14} uses automorphisms and Griffiths residues to exhibit cubic fourfolds with maximal rank of algebraic classes. In a related direction, we use finite-group actions on the Jacobian ring of a hyperplane section of $\Gr(3,10)$ to establish maximal Picard rank for Debarre--Voisin fourfolds.  We outline here the correspondence between the cohomology groups of these two spaces, which will allow us to determine the Picard ranks by computing group actions on the Jacobian ring.

Let $W_{10}$ denote a complex vector space of dimension 10, and let $\sigma \in \wedge^3 W_{10}^\vee$, $\sigma \neq 0$. To $\sigma$ we associate the varieties
\[Z_\sigma = \{W_3\in\Gr(3,W_{10}) : \sigma|_{W_3}=0\}, \qquad Y_\sigma = \{W_6\in\Gr(6,W_{10}) : \sigma|_{W_6}=0\}.\]
We call $Z_\sigma$ the \emph{Pl\"ucker hyperplane section} of $\Gr(3,10)$, and $Y_\sigma$ its associated \emph{Debarre--Voisin fourfold}. For generic $\sigma$, $Y_{\sigma}$ is a smooth irreducible hyper-K\"ahler fourfold which is of $\mathrm K3^{[2]}$ type, i.e. deformation equivalent to the Hilbert scheme of two points on a K3 surface. Debarre and Voisin \cite{Deb10} showed that these fourfolds form a 20-dimensional family. Throughout, we consider forms for which $Z_\sigma$ and $Y_\sigma$ are smooth of dimension 20 and 4, respectively.

The incidence variety

\[
    M_\sigma =\{(W_3,W_6)\in Z_\sigma\times Y_\sigma : W_3\subset W_6\}
\]
has projections
\[
    Y_\sigma\xleftarrow{q}M_\sigma\xrightarrow{p}Z_\sigma.
\]
Since the fibers of $q$ are $\Gr(3,W_6)$ of dimension 9, the correspondence $q_* p^*$ lowers the Hodge bidegree by $(9,9)$. Its restriction to vanishing cohomology induces an isomorphism of rational Hodge structures
\begin{equation}
    \label{eq: qp hodge correspondence}
    (q_*p^*)_{\van}: H^{20}_{\van}(Z_\sigma,\mathbb Q) \xrightarrow{\sim}
    H^2_{\prim}(Y_\sigma,\mathbb Q)(-9),
\end{equation}
where primitivity on $Y_\sigma$ is taken with respect to the Pl\"ucker polarization; see \cite[Corollary 2.7]{Deb10} and
\cite[Theorem 4.6]{BS21}. 

\begin{figure}[H]
    \[
        \begin{array}{ccccccccc}
            1 && 21 && 232 && 21 && 1 \\
            & 0 && 0 && 0 && 0 & \\
            && 1 && 21 && 1 && \\
            &&& 0 && 0 &&& \\
            &&&& 1 &&&&
        \end{array}
    \]
    \caption{Bottom half of the Hodge diamond for $Y_\sigma$, a hyper-K\"ahler fourfold of $K3^{[2]}$ type. Note that $h^{1,1}_{\prim}(Y_\sigma) = 20$.}
    \label{fig:hyperkahler hodge diamond}
\end{figure}

We now provide a framework for finding hypersurfaces in $\Gr(k,n)$ admitting actions of finite groups, and calculating the induced action on their vanishing cohomology. The first step is done by finding invariant lines in suitable spaces of polynomials to supply candidate hypersurfaces. After checking smoothness, we compute the induced group action on its Jacobian ring and use the residue isomorphism to interpret the result.
To construct hypersurfaces with prescribed symmetries, let $V$ be an $n$-dimensional representation of a finite group $G$. The induced action on $X= \Gr(k,V)$ also gives an action on the coordinate ring $R$ of $X$, which in degree $d$ is
\[R_{d}=H^0(X,\mathcal O_X(d))\cong\mathrm{Sym}^d({\textstyle\bigwedge^k}V^\vee)/I_d
\]
where $I_d$ is the degree $d$ part of the Pl\"ucker ideal.  A nonzero section $f\in R_d$ defines a $G$-invariant hypersurface precisely when the line $\C f$ is $G$-invariant, so finding one-dimensional subrepresentations of $R_d$ provide candidate hypersurfaces, although we must still check smoothness.

To check smoothness by reduction modulo a prime, we use the following standard fact.

\begin{lemma}[Smoothness by reduction modulo a prime]
    \label{lem: smooth reduction}
    Let $\mathcal X\to\operatorname{Spec}\Z_{(p)}$ be a proper, flat morphism of finite presentation, where $\Z_{(p)}$ is the localization of $\Z$ at a prime $p$. If its special fiber is smooth over $\mathbb F_p$, then its generic fiber is smooth over $\Q$, and hence remains smooth after base change to $\C$.
\end{lemma}

For our hypersurfaces in a Grassmannian over $\Z_{(p)}$, flatness holds whenever the defining section remains nonzero modulo $p$: its reduction is then a nonzerodivisor on the integral special fiber of the Grassmannian. We check smoothness of the hypersurface modulo $p$ by verifying that its equation and first derivatives generate the unit ideal on every standard affine chart.

For the Debarre--Voisin construction, we may take $R_1 = \bigwedge^3 W_{10}^\vee$, where $W_{10}$ is a 10-dimensional representation of a group $G$.  A $G$-invariant line $\C\sigma$ gives a pair of varieties $Z_\sigma$ and $Y_\sigma$ as before, now equipped with compatible $G$-actions which lift to the incidence variety $M_\sigma$.  Thus, the isomorphism $(q_*p^*)_\mathrm{van}$ in \eqref{eq: qp hodge correspondence} is $G$-equivariant.

By \eqref{eq: qp hodge correspondence}, the only nonzero Hodge summands of \(H^{20}_{\mathrm{van}}(Z_\sigma,\mathbb C)\) are
$$
H^{11,9}_{\mathrm{van}}(Z_\sigma)\oplus H^{10,10}_{\mathrm{van}}(Z_\sigma)\oplus H^{9,11}_{\mathrm{van}}(Z_\sigma),
$$
corresponding to the types \((2,0)\), \((1,1)\), and \((0,2)\) in \(H^2_{\mathrm{prim}}(Y_\sigma)\), shifted by \((9,9)\).  The following proposition relates the induced action on the cohomology of $Z_\sigma$ with the Picard rank of $Y_\sigma$.

\begin{proposition}
    \label{prop: maximal Picard criterion}
    Suppose a finite group $G$ acts linearly on $W_{10}$ and preserves $\mathbb C\sigma$. If $G$ acts trivially on $H^{11,9}_{\mathrm{van}}(Z_\sigma)$ and $H^{9,11}_{\mathrm{van}}(Z_\sigma)$, and
    \[
        H^{10,10}_{\mathrm{van}}(Z_\sigma)^G=0,
    \]
    then $\rho(Y_\sigma)=21$.
\end{proposition}

\begin{proof}
    The averaging idempotent
    \[
        e_G=\frac{1}{|G|}\sum_{g\in G}g\in\Q[G]
    \]
    acts on the vanishing cohomology $H^{20}_{\mathrm{van}}(Z_\sigma,\mathbb Q)$.  Since $G$ fixes $H_\mathrm{van}^{11,9}$ and $H_\mathrm{van}^{9,11}$,
    \[
        (\ker e_G)\otimes_{\mathbb Q}\mathbb C
        =
        H^{10,10}_{\mathrm{van}}(Z_\sigma).
    \]
    Thus this Hodge summand is defined over $\mathbb Q$. The isomorphism \eqref{eq: qp hodge correspondence} identifies $\ker e_G$ with a $20$-dimensional rational subspace of $H^{1,1}_{\mathrm{prim}}(Y_\sigma)$. By the Lefschetz $(1,1)$ theorem, these classes, together with the Pl\"ucker class, give Picard rank $21$.
\end{proof}

For a nonabelian simple group $G$, every one-dimensional $G$-representation is trivial. Consequently, for the groups considered below, it suffices to only check that $H_\mathrm{van}^{10,10}(Z_\sigma)$ contains no trivial summand.

The examples below illustrate how the equivariant residue description determines Hodge representations from the defining equation. In the $\mathrm{PSL}_2(11)$ case, this recovers Song's Picard rank and representation decomposition by a different method; the same computation applies to the $A_7$-invariant example.

The Picard-rank conclusion for these examples also follows from the classification of H\"ohn and Mason \cite[Theorem 8.7 and Table 9]{HM19}, once their faithful symplectic group actions are established.

Song \cite{Son21} used the group $G=\mathrm{PSL}(2,\mathbb F_{11})$ to construct a Debarre--Voisin fourfold of maximal Picard rank. For the 10-dimensional $G$-representation $V_{10}$ (as given in \cref{tab:character table PSL(2 11)}), the space of invariant alternating 3-forms is one-dimensional, generated in Pl\"ucker coordinates in a suitable basis by
\[
    \begin{aligned}
        \sigma={}& p_{025}+p_{136}+p_{247}-p_{038}-p_{149}\\
        &-p_{079}+p_{158}+p_{269}-p_{357}-p_{468}.
    \end{aligned}
\]
His proof uses the associated Peskine variety
\[
    X_\sigma = \{[v]\in\mathbb P(V_{10}) : \operatorname{rank}(v\iprod\sigma)\leq6\}.
\]
Its $55$ isolated singular points provide divisors on $Y_\sigma$, and a calculation of their Beauville--Bogomolov pairings establishes Picard rank $21$, as well as a decomposition of the primitive part of the Picard group into irreducible $G$-representations.

Below, we recover this Picard rank and $G$-representation decomposition by computing the $G$-representation on the Jacobian ring and applying \Cref{prop: maximal Picard criterion}. We then apply the same method to a hyperplane section invariant under $A_7$.

\subsection[A Debarre--Voisin Fourfold with PSL(2,F11) Symmetry]
{A Debarre--Voisin Fourfold with
\texorpdfstring{$\mathrm{PSL}(2,\mathbb F_{11})$}{PSL(2,F11)}
Symmetry}
\label{sec: our construction PSL211}

Let $G=\mathrm{PSL}(2,\mathbb{F}_{11})$, and let $V_{10}$ (as in \cref{tab:character table PSL(2 11)}) denote the irreducible representation for which
\[
    \dim\left(\bigwedge^3V_{10}^{\vee}\right)^G=1.
\]
In the coordinates obtained from the first author's Grassmannian programs \cite{GitHub:thesis}, the invariant line is generated by
\begin{equation}
    \label{eq: psl211 hyperplane}
    \begin{aligned}
        f={}&p_{017}-p_{046}+p_{059}+p_{128}-p_{156}\\
        &+p_{239}-p_{267}+p_{345}-p_{378}-p_{489}.
    \end{aligned}
\end{equation}
Since $G$ is nonabelian and simple, it fixes $\C f$ as well as the one-dimensional space $H^0(X,K_X(10))=\C\Omega$.  The permutation of coordinate covectors
\[
    x_i\longmapsto x_{\pi(i)}, \qquad
    (\pi(0),\ldots,\pi(9)) =(0,3,1,4,2,7,5,8,6,9),
\]
sends $f$ to the negative of Song's form displayed above. Consequently, $Z_f$ and its associated fourfold $Y_f$ are isomorphic to Song's examples, and are smooth \cite{Son21}.

\begin{table}[ht]
    \centering
    \renewcommand{\arraystretch}{1.15}
    \setlength{\tabcolsep}{7pt}
    \begin{adjustbox}{max width=\textwidth}
        \setlength{\tabcolsep}{7pt}
        \begin{tabular}{|c| r r r r r r r r|}
            \hline
            &   $(1)$ & $b^{-1}aba^{-1}b^{-1}a^{-2}$ & $b^{-1}a^{-1}$ & $aba^{-1}b^{-1}a^{-2}b$ & $b^{-1}$ & $a^{-2}b^{-1}a^{-2}$ & $a$ & $a^2$\\
            \hline
            Size & 1 & 60 & 60 & 55 & 110 & 110 & 132 & 132\\
            \hline
            $\mathbbm{1}$ & $1$  & $1$ & $1$ & $1$ & $1$ & $1$ & $1$ & $1$ \\
            $V_5$ & $5$  & $\frac{-1+\sqrt{-11}}{2}$ & $\frac{-1-\sqrt{-11}}{2}$ & $1$ & $-1$ & $1$ & $0$ & $0$ \\
            $V_5^\vee$ & $5$  & $\frac{-1-\sqrt{-11}}{2}$ & $\frac{-1+\sqrt{-11}}{2}$ & $1$ & $-1$ & $1$ & $0$ & $0$ \\
            $V_{10}$ & $10$ & $-1$ & $-1$ & $-2$ & $1$ & $1$ & $0$ & $0$ \\
            $V_{10}'$ & $10$ & $-1$ & $-1$ & $2$ & $1$ & $-1$ & $0$ & $0$ \\
            $V_{11}$ & $11$ & $0$ & $0$ & $-1$ & $-1$ & $-1$ & $1$ & $1$ \\
            $V_{12}$ & $12$ & $1$ & $1$ & $0$ & $0$ & $0$ & $\frac{-1-\sqrt{5}}{2}$ & $\frac{-1+\sqrt{5}}{2}$ \\
            $V_{12}'$ & $12$ & $1$ & $1$ & $0$ & $0$ & $0$ & $\frac{-1+\sqrt{5}}{2}$ & $\frac{-1-\sqrt{5}}{2}$ \\
            \hline
        \end{tabular}
    \end{adjustbox}
    \caption{Character table of $\mathrm{PSL}(2,\mathbb F_{11})$. Representatives are expressed in terms of $a=(3\,11\,9\,7\,5)(4\,12\,10\,8\,6)$ and $b=(1\,2\,8)(3\,7\,9)(4\,10\,5)(6\,12\,11)$. Products of permutations act from left to right.}
    \label{tab:character table PSL(2 11)}
\end{table}

Let $X = \Gr(3,V_{10})$, hence $\dim X = 21$. Let $R/J_f$ denote the Jacobian ring of $X$ associated to $f$. The hypotheses of \Cref{thm: main theorem} hold because $\gcd(3,10) = 1$ and $3\notin\{2,8\}$. Moreover, $H^{21}_{\prim}(X)=0$ since a complex Grassmannian has no odd cohomology. Hence the residue isomorphism gives
\[
    (R/J_f)_0\cong H^{11,9}_{\mathrm{van}}(Z_f), \qquad
    (R/J_f)_1\cong H^{10,10}_{\mathrm{van}}(Z_f),\qquad
    (R/J_f)_2\cong H^{9,11}_{\mathrm{van}}(Z_f),
\]
of dimensions 1,20, and 1, respectively.  These identifications are $G$-equivariant.

We now investigate the action on $(R/J_f)_1$.  Under the infinitesimal action map $\operatorname{End}(V_{10})\longrightarrow R_1$ given by $A\longmapsto\zeta_A(f)$, the scalar matrices have image $\mathbb C f$, which together with the images of trace-free elements spans $(J_f)_1$.  Applying the matrix units to $f$ in \texttt{debarre\_voisin.sage} gives a matrix of rank $100$ modulo $113$, so the action map is injective in characteristic zero. Thus,
\[
    \chi_{(R/J_f)_1}=\chi_{\bigwedge^3 V_{10}^{\vee}}-\chi_{\operatorname{End}(V_{10})}.
\]
Evaluating on the conjugacy classes gives the character
\[
    \chi_{(R/J_f)_1} =(20,-2,-2,4,2,-2,0,0),
\]
with columns ordered as in \Cref{tab:character table PSL(2 11)}. This is twice the character of $V_{10}'$, so
\[
    H^{10,10}_{\mathrm{van}}(Z_f)\cong (R/J_f)_1\cong (V_{10}')^{\oplus2}.
\]
$G$ acts trivially on the two outer Hodge summands, so \Cref{prop: maximal Picard criterion} gives
$$\rho(Y_f) = 21.
$$

\subsection[A Debarre--Voisin Fourfold with A7 Symmetry]
{A Debarre--Voisin Fourfold with
\texorpdfstring{$A_7$}{A7} Symmetry}
\label{sec: a7 sym}
We now apply the methods of the previous section to the alternating group $G=A_7$, which has the character table given in \Cref{tab:A7 character table}. Other fourfolds with $A_7$ symmetry were constructed by Billi and Wawak \cite[Theorem 0.1]{BW24} as double EPW sextics; here we consider a Debarre--Voisin model. 

\begin{table}
    \[
        \newcommand{\bseven}{\frac{-1+\sqrt{-7}}{2}}
        \newcommand{\bsevenbar}{\frac{-1-\sqrt{-7}}{2}}
        \begin{array}{|c|rrrrrrrrr|}
            \hline
            & 1 & ab^{-1}ab & a & a^{-1}bab & a^{-1}bab^2 & b & ababab^2 & ab & a^{-1}b \\
            \hline
            \text{Size} & 1 & 105 & 70 & 280 & 630 & 504 & 210 & 360 & 360 \\
            \hline
            \mathbbm{1} & 1  &  1 &  1 &  1 &  1 &  1 &  1 &  1 &  1 \\
            V_6         & 6  &  2 &  3 &  0 &  0 &  1 & -1 & -1 & -1 \\
            V_{10}      & 10 & -2 &  1 &  1 &  0 &  0 &  1 & \bseven & \bsevenbar \\
            V_{10}'     & 10 & -2 &  1 &  1 &  0 &  0 &  1 & \bsevenbar & \bseven \\
            V_{14}      & 14 &  2 &  2 & -1 &  0 & -1 &  2 &  0 &  0 \\
            V_{14}'     & 14 &  2 & -1 &  2 &  0 & -1 & -1 &  0 &  0 \\
            V_{15}      & 15 & -1 &  3 &  0 & -1 &  0 & -1 &  1 &  1 \\
            V_{21}      & 21 &  1 & -3 &  0 & -1 &  1 &  1 &  0 &  0 \\
            V_{35}      & 35 & -1 & -1 & -1 &  1 &  0 & -1 &  0 &  0 \\
            \hline
        \end{array}
    \]
    \caption{Character table for $A_7$. Representatives of each conjugacy class are given in terms of the generators $a = (1\,2\,3)$ and $b = (3\,4\,5\,6\,7)$.}
    \label{tab:A7 character table}
\end{table}

$A_7$ has two complex-conjugate irreducible representations of dimension 10, each satisfying
\[
    \dim\left(\bigwedge^3 V^\vee\right)^{G}=1.
\]
We will use the representation $V_{10}$. In a suitable basis, the $A_7$-invariant line in $R_1$ is given by
\begin{equation}
    \label{eq: a7 hyperplane}
    \begin{aligned}
        f={}&p_{017}+p_{028}+p_{039}+9p_{123}\\
        &-5p_{126}+5p_{135}-5p_{234}\\
        &+2p_{149}+2p_{257}+2p_{368}
        +4p_{456}+p_{789}.
    \end{aligned}
\end{equation}

Let $Z_f\subset\Gr(3,V_{10})$ be the corresponding hypersurface.  The script \texttt{debarre\_voisin.sage} checks smoothness modulo 113, so \Cref{lem: smooth reduction} gives smoothness over $\Q$, and hence over $\C$. The associated Debarre--Voisin variety $Y_f$ is therefore a smooth fourfold \cite[Lemma 3.1]{BS21}.

As in \Cref{sec: our construction PSL211}, the residue isomorphism identifies 

\[
    (R/J_f)_0\cong H^{11,9}_{\mathrm{van}}(Z_f), \qquad
    (R/J_f)_1\cong H^{10,10}_{\mathrm{van}}(Z_f),\qquad
    (R/J_f)_2\cong H^{9,11}_{\mathrm{van}}(Z_f),
\]
of dimensions 1,20, and 1, respectively. These identifications are $A_7$-equivariant: both $f$ and the twisted volume form $\Omega$ are fixed, since $A_7$ has no nontrivial one-dimensional representations.

The derivative matrix again has rank $100$ modulo $113$. Subtracting the character of $\operatorname{End}(V_{10})$ from that of $\bigwedge^3V_{10}^{\vee}$, as above, gives
\[
    \chi_{(R/J_f)_1}
    =(20,4,2,2,0,0,-2,-1,-1),
\]
with columns ordered as in
\Cref{tab:A7 character table}. This is the sum of the characters of $V_6$ and $V_{14}'$, so
\[
    H^{10,10}_{\mathrm{van}}(Z_f) \cong (R/J_f)_1 \cong V_6\oplus V_{14}'.
\]
Neither summand is trivial, so by \Cref{prop: maximal Picard criterion},
\[\rho(Y_f)=21.
\]

\subsection{An Explicit Picard--Fuchs Equation for a Symmetric Gushel--Mukai Fivefold}
We now use Griffiths--Dwork reduction to compute a Picard--Fuchs equation. We consider a pencil of symmetric Gushel--Mukai fivefolds, whose Jacobian rings are described in \cite{Fat21}. Let $V_5$ be the irreducible five-dimensional rational representation of $G=A_5$, and set $X=\Gr(2,V_5)$.  The space of invariant quadrics on $X$ has dimension three.  We choose two linearly independent invariant quadrics to form our pencil, and use symmetry to restrict the period calculation to a two-dimensional subspace of the middle cohomology.

\begin{table}[ht]
    \centering
    \renewcommand{\arraystretch}{1.15}
    \setlength{\tabcolsep}{7pt}
    \begin{adjustbox}{max width=\textwidth}
        \begin{tabular}{|c|rrrrr|}
            \hline
            & $1$ & $a$ & $b$ & $ab$ & $(ab)^2$ \\
            \hline
            Size & $1$ & $15$ & $20$ & $12$ & $12$ \\
            \hline
            $\mathbbm{1}$ & $1$ & $1$ & $1$ & $1$ & $1$ \\
            $V_3$ & $3$ & $-1$ & $0$
            & $\frac{1+\sqrt5}{2}$ & $\frac{1-\sqrt5}{2}$ \\
            $V_3'$ & $3$ & $-1$ & $0$
            & $\frac{1-\sqrt5}{2}$ & $\frac{1+\sqrt5}{2}$ \\
            $V_4$ & $4$ & $0$ & $1$ & $-1$ & $-1$ \\
            $V_5$ & $5$ & $1$ & $-1$ & $0$ & $0$ \\
            \hline
        \end{tabular}
    \end{adjustbox}
    \caption{Character table for $A_5$. Representatives are
    expressed in terms of the generators $a=(1\,2)(3\,4)$
    and $b=(1\,3\,5)$.}
    \label{tab:A5 character table}
\end{table}

In the coordinates constructed by the accompanying script \texttt{pencil\_calculations.sage}, two invariant quadrics are
\[
    \begin{aligned}
        f_0={}&5p_{12}^2+(p_{14}+p_{24}+p_{35})^2
        +(p_{15}+p_{34})^2+(p_{23}+p_{45})^2,\\[4pt]
        f_\infty={}&-9p_{12}^2+(2p_{13}+p_{23})^2+(p_{14}-p_{24})^2
        +(p_{15}+2p_{25})^2\\
        &\qquad-p_{34}^2-p_{35}^2-p_{45}^2
        -6\bigl(p_{15}p_{34}+(p_{14}+p_{24})p_{35}+p_{23}p_{45}\bigr).
    \end{aligned}
\]
These generators define the fibers at $[1:0]$ and $[0:1]$ of the pencil over $\mathbb P^1_{\mathbb C}$ given by
\[
    Z_{[u:v]}=\{u f_0+v f_\infty=0\}\subset X.
\]
On the affine chart $u\ne0$, we write $t=v/u$ and denote
\[
    f_t=f_0+t f_\infty,\qquad Z_t=\{f_t=0\}\subset X.
\]
We chose these quadrics for their relative simplicity, so that the reduction calculations could be demonstrated explicitly.  As the fibers at $t=0$, $t=1$ and $t=\infty$ are singular, we check smoothness in the convenient case $t=2$.  The accompanying script \texttt{pencil\_calculations.sage} checks smoothness on all ten standard affine charts of $Z_2$ over $\mathbb F_{101}$, so \Cref{lem: smooth reduction} gives smoothness in characteristic zero.  Now the subset $B\subset\mathbb P^1_{\mathbb C}$ parameterizing smooth members of the pencil forms a nonempty Zariski-open subset, so the pencil restricts to a smooth family of Gushel--Mukai fivefolds over $B$.  We shall continue to use $t=2$ when a smooth specialization is needed.

For $t\in B$, we apply \Cref{thm: main theorem} with $(k,n,d)=(2,5,2)$ and $N=6$, which satisfies the arithmetic hypotheses.  The middle cohomology of $Z_t$ is entirely vanishing cohomology, since $H^7(X)=0$. At pole order $p=3$, the cohomology classes are represented by forms $g\Omega/f_t^p$ where $g\in R_{pd-n}=R_1$, so the residue isomorphism gives
\begin{equation}\label{eq: 3,2 cohomology of Z_t}
    H^{3,2}(Z_t)\cong(R/J_{f_t})_1=R_1\cong\bigwedge\nolimits^2 V_5^\vee\cong V_3\oplus V_3'\oplus V_4.
\end{equation}
There is no ambient contribution since $H^6_{\mathrm{prim}}(X)=0$, and the decomposition at the end is a quick character calculation.  Since the group fixes $f_t$ and $\Omega$, the residue isomorphism is $A_5$-equivariant.

At pole orders $p=1,2$, \Cref{thm: main theorem} gives $H^{5,0}(Z_t)=H^{4,1}(Z_t)=0$ (the Jacobian ring is zero in the corresponding degrees).  Thus $h^{3,2}(Z_t)=h^{2,3}(Z_t)=10$.  The group action commutes with complex conjugation, so as $G$-representation $H^{2,3}(Z_t)$ is the conjugate of $H^{3,2}(Z_t)$. Since the characters of $V_3,V_3'$ and $V_4$ are real, both Hodge summands have the same decomposition as $G$-representations.  The middle row of the Hodge diamond of $Z_t$ is given in \Cref{fig:gm fivefold hodge diamond}.
\begin{figure}[H]
    \[
        \begin{array}{ccccccccccc}
            0 && 0 && 10 && 10 && 0 && 0 \\
            & 0 && 0 && 2 && 0 && 0 & \\
            && 0 && 0 && 0 && 0 && \\
            &&& 0 && 1 && 0 &&& \\
            &&&& 0 && 0 &&&& \\
            &&&&& 1 &&&&&
        \end{array}
    \]
    \caption{Bottom half of the Hodge diamond for $Z_t$, a smooth Gushel--Mukai fivefold.}
    \label{fig:gm fivefold hodge diamond}
\end{figure}

Let $A_4\subset A_5$ be the standard inclusion.  Restriction of characters gives
\[
    \dim V_3^{A_4}=\dim(V_3')^{A_4}=0,\qquad \dim V_4^{A_4}=1.
\]
Thus $H^{3,2}(Z_t)$ and $H^{2,3}(Z_t)$ both have a one-dimensional $A_4$-invariant subspace, and we have $\dim(H^5(Z_t,\mathbb C)^{A_4})=2$.  Under the isomorphism \eqref{eq: 3,2 cohomology of Z_t}, the subspace $V_4^{A_4}$ corresponds to the one-dimensional space
\begin{equation}\label{eq: invariant line R1}
    R_1^{A_4}=\C p_{12}.
\end{equation}
We take
\[
    \omega_t=\Res_{Z_t}\frac{2!\,p_{12}\Omega}{f_t^3}\in H^{3,2}(Z_t).
\]
Let $\nabla_t=\nabla_{\partial/\partial t}$ denote the derivative in the $t$-direction. Since the numerator of $\omega_t$ is independent of $t$, differentiation gives
\[
    \nabla_t\omega_t=\Res_{Z_t}\frac{-3!\,p_{12}f_\infty\Omega}{f_t^4},\qquad
    \nabla_t^2\omega_t=\Res_{Z_t}\frac{4!\,p_{12}f_\infty^2\Omega}{f_t^5}.
\]
Since $\dim(H^5(Z_t,\mathbb C)^{A_4})=2$, the forms $\omega_t$, $\nabla_t\omega_t$, and $\nabla_t^2\omega_t$ are linearly dependent.  We apply Griffiths--Dwork reduction to produce a relation between them.

Fix a basis $\{A_i\}$ of $\mathfrak{sl}(V_5)$, and use the notation $[g]_p$ of \Cref{def: normalized pole class} with $f=f_t$.  All polynomial identities below are taken in the Pl\"ucker coordinate ring over $\C(t)$.  Since $A_4$ fixes $f_t$ and $\Omega$, \Cref{cor: invariant reduction} allows us to perform reduction without leaving the invariant subspaces.  Averaging the standard monomials over $A_4$ in \texttt{pencil\_calculations.sage} gives
\[
    \dim R_1^{A_4}=1,\qquad \dim R_3^{A_4}=15,\qquad \dim R_5^{A_4}=96.
\]

We choose complements $W_m$ of $(J_{f_t})_m^{A_4}$ in $R_m^{A_4}$ for $m=5,3$.  But $H^{1,4}(Z_t)=0$ by \Cref{fig:gm fivefold hodge diamond}, so $(J_{f_t})_5^{A_4}=R_5^{A_4}$ and hence $W_5=0$.  In \texttt{pencil\_calculations.sage}, it is shown that $p_{12}f_\infty$ is nonzero in $(R/J_{f_t})_3$ for $t=2$, hence for generic $t$.  Since $(R/J_{f_t})_3^{A_4}\cong H^{2,3}(Z_t)^{A_4}$ is one-dimensional, we take $W_3=\mathrm{Span}(p_{12}f_\infty)$.

Start with $g_5=p_{12}f_\infty^2$.  Since $W_5=0$, we set $a_5=0$ and $b_5=g_5$.  We write $b_5$ in terms of the generators of $J_{f_t}$ and apply \eqref{eq: factorial reduction} to produce an invariant cubic $g_4$ with $[g_5]_5=[g_4]_4$.  Next, we decompose $g_4$ into a scalar multiple of $p_{12}f_\infty$ plus an element in $J_{f_t}$:
\[
    \begin{aligned}
        a_4&=\lambda(t)(-p_{12}f_\infty),\\
        b_4&=c_4f_t+\sum_i h_4^i\zeta_{A_i}(f_t),\qquad c_4,h_4^i\in R_1\otimes\C(t)
    \end{aligned}
\]
(the sign on $p_{12}f_\infty$ is for convenience).  Reduction gives $[g_4]_4=[a_4]_4+[g_3]_3$, where
\[
    g_3=3c_4+\sum_i\zeta_{A_i}(h_4^i).
\]
Now $g_3$ is $A_4$-invariant and linear, so \eqref{eq: invariant line R1} gives $g_3=\mu(t)p_{12}$ for some nonzero $\mu(t)\in\C(t)$.  Since $(J_{f_t})_1=0$, we have $a_3=g_3$.  Combining the reduction with $\Res_{Z_t}[-p_{12}f_\infty]_4=\nabla_t\omega_t$ and $\Res_{Z_t}[p_{12}]_3=\omega_t$ yields
\[
    \nabla_t^2\omega_t=\lambda(t)\nabla_t\omega_t+\mu(t)\omega_t.
\]
The script \texttt{pencil\_calculations.sage} calculates
\begin{align}
    \lambda(t)&=-\frac1t-\frac1{t-1}-\frac5{5t-3}-\frac9{9t-5}-\frac6{3t-1}-\frac9{9t+1},\label{eq: new pencil lambda}\\
    \mu(t)&=-\frac{3645t^4-5022t^3+1944t^2-222t-9}{t(t-1)(3t-1)^2(5t-3)(9t+1)}.\label{eq: new pencil mu}
\end{align}

At $t=2$, the class $\omega_2$ is nonzero in $F^3H^5(Z_2,\C)$, while $\left.\nabla_t\omega_t\right|_{t=2}$ has nonzero image in the associated graded $\operatorname{Gr}_F^2H^5(Z_2,\C)$.  Hence $\omega_t$ and $\nabla_t\omega_t$ are linearly independent over $\C(t)$.  Thus, we have proved the following.

\begin{proposition}\label{prop: new pencil PF}
    For the pencil $f_t=f_0+t f_\infty$, the $A_4$-invariant residue $\omega_t$ satisfies the minimal second-order Picard--Fuchs equation
    \[
        \nabla_t^2\omega_t=\lambda(t)\nabla_t\omega_t+\mu(t)\omega_t,
    \]
    with $\lambda$ and $\mu$ as in \eqref{eq: new pencil lambda}--\eqref{eq: new pencil mu}.  In particular, for any locally flat family of five-cycles $\gamma_t$, the period $\Pi(t)=\int_{\gamma_t}\omega_t$ satisfies
    \[
        \Pi''=\lambda(t)\Pi'+\mu(t)\Pi.
    \]
\end{proposition}

Although the reduced representatives depend on the complements $W_m$ chosen in the reduction algorithm, the coefficients $\lambda(t)$ and $\mu(t)$ are independent of those complements.  Indeed, for the fixed section $\omega_t$ and parameter $t$, the classes $\omega_t$ and $\nabla_t\omega_t$ form a basis of the $A_4$-invariant cohomology over $\mathbb C(t)$, so their coefficients in the expression for $\nabla_t^2\omega_t$ are unique.
\vfill

\printbibliography

\end{document}